# Обобщения теоремы Штейнера-Лемуса о признаках равнобедренности треугольника

Рабе Алексей Дмитриевич [1], Привалов Александр Андреевич [2]

[1] ГБОУ г. Москвы "Школа на Юго-Востоке имени Маршала В.И. Чуйкова"

[2] This paper is prepared under the supervision of Alexander Privalov and is submitted to the Moscow Mathematical Conference for High-School Students. Readers are invited to send their remarks and reports on this paper to mmks@mccme.ru'

## Abstract

In 1840 Jacob Steiner on Christian Rudolf's request proved that a triangle with two equal bisectors is isosceles. But what about changing the bisectors to cevians? Cevian is any line segment in a triangle joining a vertex of the triangle and a point on the opposite side. Not for any pairs of equal cevians the triangle is isosceles.

**Theorem.** *If for a triangle ABC there are equal cevians issuing from A and B, which intersect on the bisector or on the median of the angle C, then AC=BC (so the triangle ABC is isosceles).*

**Proposition.** *Let ABC be an isosceles triangle. Define circle $\omega$ to be the circle symmetric relative to AB to the circumscribed circle of the triangle ABC. Then the locus of intersection points of pairs of equal cevians is the union of the base AB, the triangle's axis of symmetry, and the circle $\omega$.*

## 1   Исторический обзор

Формулировки основных результатов находятся в параграфе 2.

Данная работа является продолжением исследований, начатых Я. Штейнером, К. Лемусом и А. Ботемой о признаках равнобедренности треугольника. Требовалось выяснить, возможно ли обобщить теорему Штейнера-Лемуса, сформулированную как один из признаков равнобедренного треугольника: *если в треугольнике две биссектрисы равны, то он является равнобедренным.*

В ходе проведённых исследований было доказано несколько признаков равнобедренности треугольника, вытекающих из равенства некоторых его чевиан (*чевианой* треугольника называют отрезок, соединяющий вершину с точкой, лежащей на противоположной стороне или ее продолжении и отличной от вершин этой стороны).

Примерами чевиан являются биссектрисы, медианы и высоты. Нетрудно показать, что треугольники с двумя равными высотами или с двумя равными медианами являются равнобедренными. Этот простой факт был, вероятно, известен ученым Древней Греции. Однако с биссектрисами дело обстояло несколько иначе.

В 1840 году Кристианом Людольфом Лемусом был задан вопрос Якобу Штейнеру о возможных доказательствах следующей теоремы, которая вошла в историю как

**Теорема Штейнера-Лемуса [3].** *Треугольник с двумя равными биссектрисами является равнобедренным.*

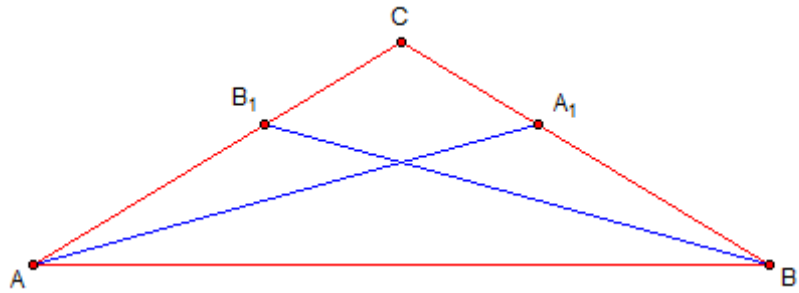

Но будет ли треугольник равнобедренным, если заменить равные биссектрисы $AA_1$ и $BB_1$ на равные чевианы? То есть, поставим следующую задачу: *выяснить, из равенства каких чевиан треугольника будет вытекать равенство его сторон*. В работе рассмотрены аналоги теоремы Штейнера-Лемуса о признаках равнобедренности треугольника по равными чевианами.

## 2    Формулировки основных результатов

Чевианой треугольника называют отрезок, соединяющий вершину с точкой, лежащей на противоположной стороне или ее продолжении и отличной от вершин этой стороны. А чтобы избежать путаницы в обозначениях прямой и отрезка, будем обозначать прямую, соединяющие точки $A$ и $B$, через $(AB)$.

**Теорема 1.** *Если в треугольнике ABC существуют равные чевианы, выходящие из A и B, пересекающиеся на биссектрисе угла C, то такой треугольник равнобедренный.*

**Теорема 2.** *Если в треугольнике ABC существуют равные чевианы, выходящие из A и B, пересекающиеся на медиане угла C, то такой треугольник равнобедренный.*

Для треугольника $ABC$ обозначим как $C_{ABC}$ окружность, симметричную окружности, описанной около треугольника $ABC$, относительно $(AB)$.

**Теорема 3.** *Если в треугольнике ABC существуют равные чевианы, выходящие из A и B, пересекающиеся на окружности $C_{ABC}$, то такой треугольник равнобедренный.*

**Теорема 4.** *Если в треугольнике ABC существуют равные чевианы, выходящие из A и B, пересекающиеся на высоте CH в точке O, причем* $y = \dfrac{\overrightarrow{HO} \cdot \overrightarrow{HC}}{|HC|}$ *не удовлетворяет уравнению*

$$y^3 + y(ctg^2\alpha + ctg^2\beta - 1) + 2ctg\alpha \cdot ctg\beta = 0$$

*где* $\alpha = \angle BAC$ *и* $\beta = \angle ABC$, *то такой треугольник равнобедренный.*

Во избежание неоднозначностей в трактовке, формулируя и доказывая следующую теорему, будем считать все углы ориентированными. Тогда для $\angle BAC$ луч $AA_1$ назовем $k$-трисой ($-\infty < k < \infty$), если $\angle BAA_1 = k\angle BAC$.

**Теорема 5 (обобщение теоремы Штейнера-Лемуса).** *Если равные чевианы $AA_1$ и $BB_1$ являются $k$-трисами треугольника ABC и $0 < k \le 1$, то треугольник ABC равнобедренный.*

**Теорема 6.** *Пусть шесть чевиан треугольника ABC равны ($AA_1= AA_2= BB_1= BB_2= CC_1= CC_2$). Тогда точки $A_1, A_2, B_1, B_2, C_1, C_2$ принадлежат одной конике.*

**Теорема 7.** *Геометрическое место точек пересечения пар равных чевиан $AA_1$ и $BB_1$ треугольника ABC есть объединение отрезка AB и кубической кривой.*

Тетраэдр называется равногранным, если все его грани – равновеликие треугольники.
**Гипотеза 8 (обобщение теоремы Штейнера-Лемуса в трехмерном случае).** *Если в тетраэдре все биссектрисы трехгранных углов равны, то такой тетраэдр равногранный.*

## 3    Доказательство теоремы 1

Для доказательства этой теоремы понадобятся следующие две леммы:

**Лемма 1.1.** *Треугольники равны, если у них равна сторона, противолежащий угол и биссектриса этого угла.*

*Доказательство леммы 1.1.* Очевидно, что все треугольники, имеющие по равной стороне и по равному, противолежащей ей углу, вписаны в одну окружность. Биссектрисы этих углов пересекаются точке $P$, лежащей на этой окружности. На рисунке 2 $CD$ и $C_1D_1$ – биссектрисы двух таких треугольников $ABC$ и $ABC_1$. Стоит указать, что точки $C$ и $C_1$ лежат по одну сторону от диаметра, перпендикулярного AB.

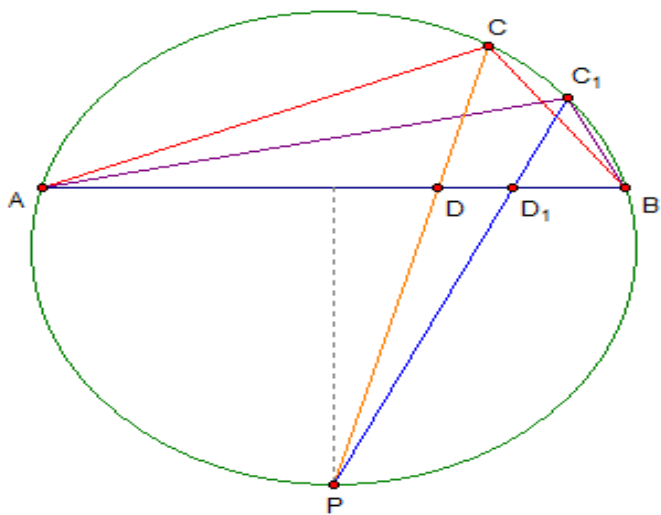

Рис. 2

Тогда следующее соотношение:

$$CD = PC - PD > PC_1 - PD > PC_1 - PD_1 = C_1D_1 \qquad (1)$$

показывает, что разным биссектрисам соответствуют разные треугольники.
  Лемма 1.1 доказана.

**Лемма 1.2.** *Треугольники равны, если у них равны сторона, противолежащий угол и биссектриса внешнего к этому углу угла.*

*Доказательство леммы 1.2.* Рассмотрим треугольники $ABC$ и $ABC_1$ такие, что $\angle ACB = \angle AC_1B$. Хорошо известно, что биссектриса внешнего угла треугольника перпендикулярна биссектрисе этого угла. Отсюда следует, что если биссектриса внешнего к $\angle ACB$ угла не пересекает сторону $(AB)$, то биссектриса угла $ACB$ будет являться высотой треугольника $ABC$. Следовательно, треугольник $ABC$ равнобедренный, как и треугольник $ABC_1$. Поэтому, треугольник $ABC$ равен треугольнику $ABC_1$ (по второму признаку равенства треугольников). Стоит указать, что точки $C$ и $C_1$ лежат по одну сторону от диаметра, перпендикулярного $AB$.

Пусть эти биссектрисы пересекают сторону $(AB)$. Обозначим их $CE$ и $C_1E_1$ (рис.3).

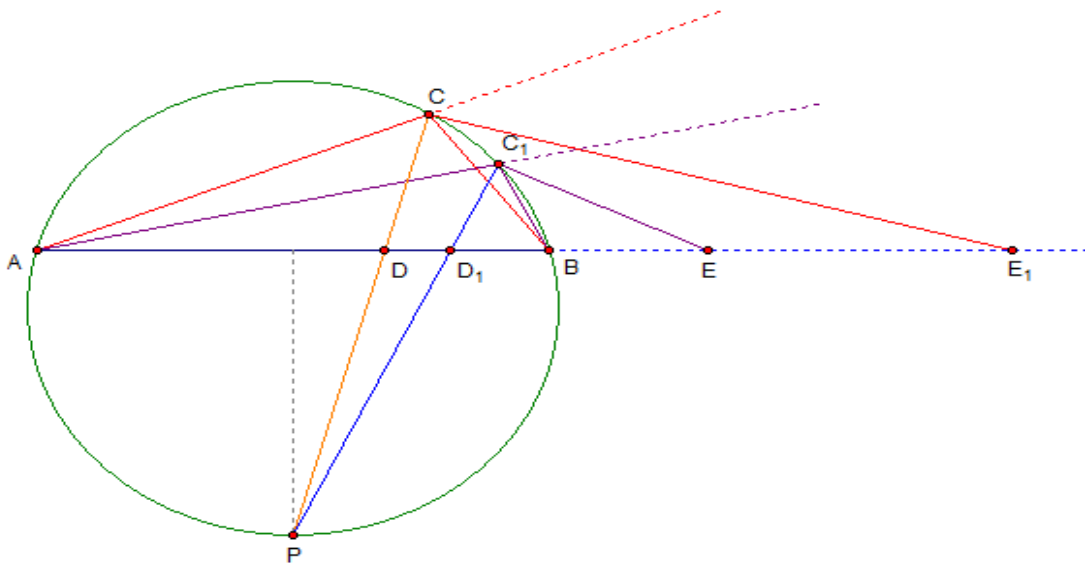

Рис. 3

Из (1) в лемме 1 следует, что катет $CD$ прямоугольного треугольника $DCE$ больше катета $C_1D_1$ треугольника $D_1C_1E_1$, кроме того, $\angle CDE > \angle C_1D_1E_1$, следовательно $CE > C_1E_1$, т.е. неравным биссектрисам соответствуют неравные треугольники.

Лемма 1.2 доказана.

*Доказательство теоремы 1.* Рассмотрим первый случай: отрезки $AA_1$ и $BB_1$ параллельны сторонам $CB$ и $CA$ соответственно и пересекаются в точке $O$, лежащей на биссектрисе (рис.4). Тогда четырехугольник $OABC$ является параллелограммом с диагональю $CO$, являющейся биссектрисой. Следовательно, $OABC$ – ромб и $CA=CB$.

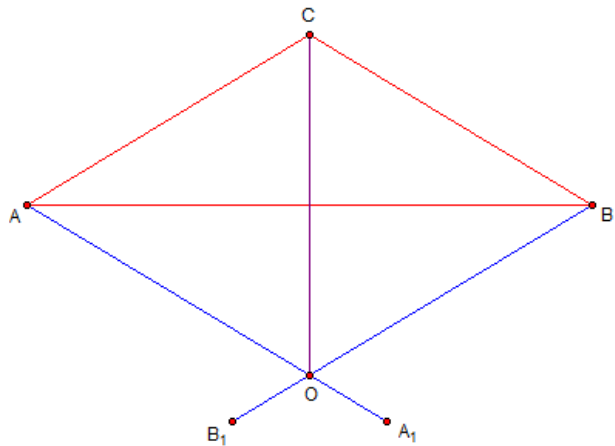

Рис. 4

Второй случай: чевианы $AA_1$ и $BB_1$ треугольника $ABC$ равны и пересекаются на биссектрисе ($CD$) в точке $O$ так, что точки $A_1$ и $B_1$ лежат по разные стороны с точкой $C$ относительно $(AB)$ (рис.5).

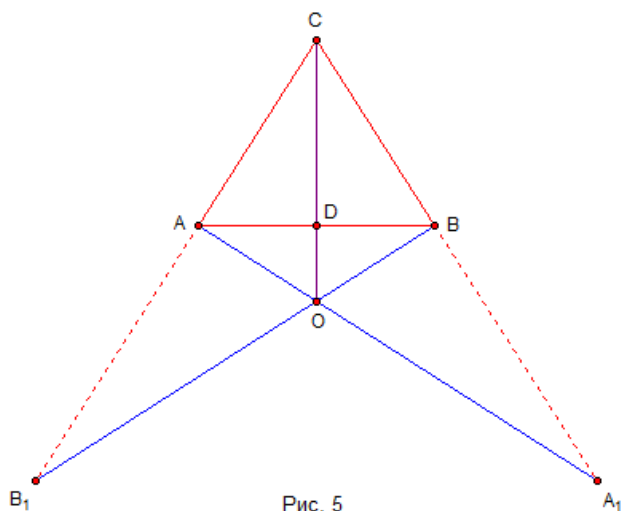

Рис. 5

Здесь, как и выше, из леммы 1.1 следует равенство треугольников $ACA_1$ и $BCB_1$, т.е. $CA=CB$ и треугольник $ABC$ – равнобедренный.

Третий случай: чевианы $AA_1$ и $BB_1$ треугольника $ABC$ равны и пересекаются на биссектрисе ($CD$) в точке $O$ так, что точки $A_1$ и $B_1$ лежат по одну сторону с точкой $C$ относительно $(AB)$. Для доказательства третьего случая рассмотрим рисунки 6а и 6б.

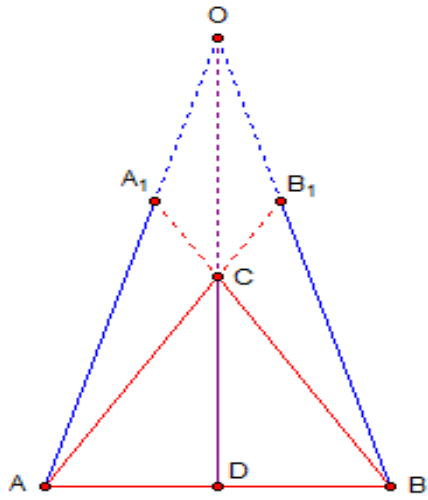
Рис. 6а

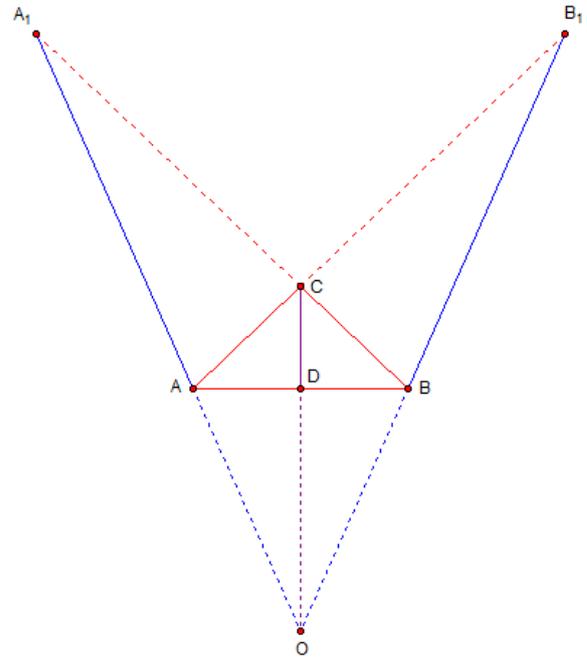
Рис. 6б

Здесь треугольники $ACA_1$ и $BCB_1$ имеют равные стороны ($AA_1=BB_1$), равные противолежащие углы ($\angle ACA_1=\angle BCB_1$) общую равную биссектрису $CO$ внешних углов $\angle ACB$. По лемме 1.2 эти треугольники равны, т.е. $AC = BC$.

Теорема 1 доказана.

**Замечание 1.3.** Ботема рассматривал теорему Штейнера-Лемуса для биссектрис внешних углов треугольника. Он построил треугольник $ABC$, у которого биссектрисы $AA_1$ и $BB_1$ внешних углов равны и равны стороне $AB$. В этом случае точки $A_1$ и $B_1$ лежат по разные стороны прямой ($AB$).

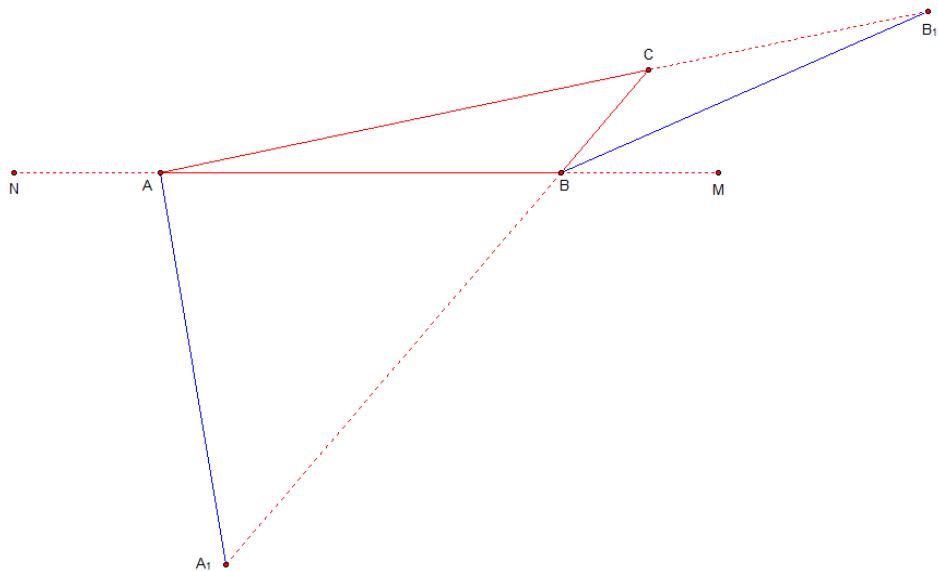

*Доказательство.* Обозначим угол $\angle CAB=\alpha$. Тогда, т.к. $AB=BB_1$, то $\angle CB_1B=\alpha$. Отсюда и того, что $BB_1$ – биссектриса угла $CBM$, то $\angle B_1BM=2\alpha$ и $\angle CBM =4\alpha$. Так как $AB=AA_1$, то $\angle AA_1B=\angle ABA_1=\angle CBM=4\alpha$ и $\angle BAA_1=\pi-8\alpha$. Далее, $AA_1$ – биссектриса угла $\angle BAN = \pi-\alpha$, отсюда имеем $\pi-\alpha=2(\pi-8\alpha)$ и $\angle CAB=\alpha=12^0$, $\angle ACB=3\alpha=36^0$ и $\angle ABC=132^0$.

# 4 Доказательство теоремы 2

*Первое доказательство теоремы 2.* Если чевианы $AA_1$ и $BB_1$ пересекаются на медиане, то отрезок $A_1B_1$ параллелен основанию треугольника $AB$. Тогда трапеция $ABA_1B_1$ имеет равные диагонали и, следовательно, является равнобокой (рис. 8).

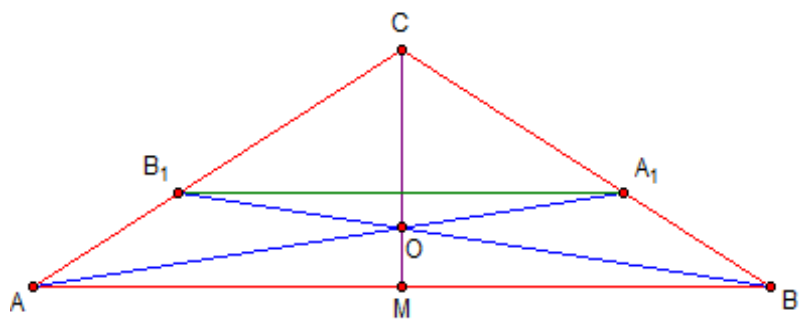

Рис. 8

*Второе доказательство теоремы 2.* Пусть точка $O$ лежит на медиане $CM$ треугольника $ABC$. Тогда ее барицентрические координаты будут $O(1 : 1 : x)$. Несложно видеть, что если $x = -1$, то $AOBC$ – параллелограмм; если $x = 0$, то точка $O$ совпадет с точкой $M$. Эти значения $x$ исключены в связи с тривиальностью решения.

Сначала рассмотрим случай, когда точка $O$ бесконечно удалена. Тогда отрезки $AA_1$ и $BB_1$ параллельны и по условию равны. Отсюда следует, что треугольники $ACA_1$ и $BCB_1$ имеют равную сторону и равные углы. Следовательно, они равны и треугольник $ABC$ равнобедренный (рис. 9).

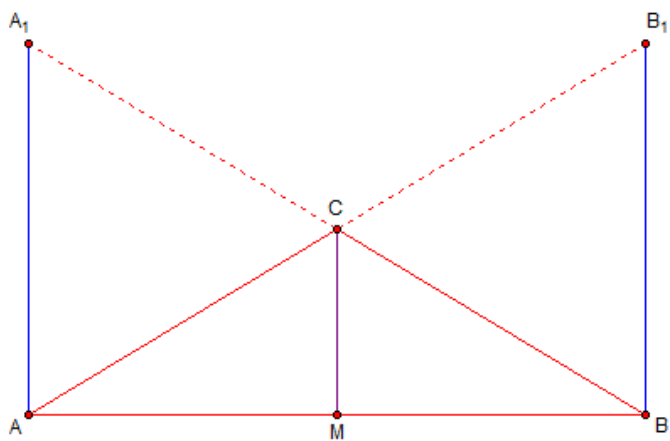

Рис. 9

Пусть теперь точка $O$ лежит на $(CM)$ и является центром масс системы материальных точек $\{1A, 1B, xC\}$, тогда точка $A_1$ является барицентром системы материальных точек $\{1B, xC\}$, а $B_1$ –центр масс системы $\{1A, xC\}$ (рис. 10).

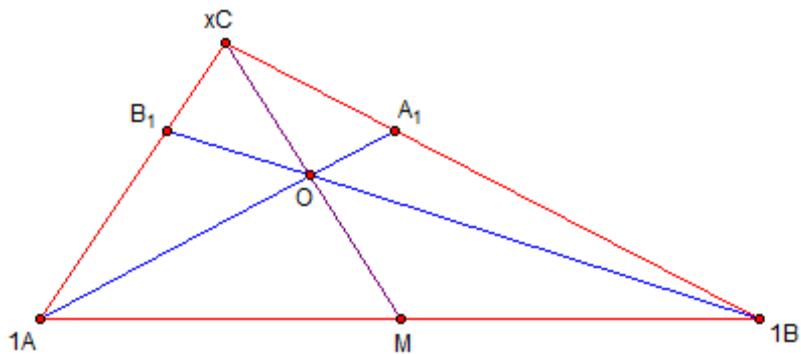

Рис. 10

Отсюда найдем векторы $\overrightarrow{AA_1}$ и $\overrightarrow{BB_1}$:

$$\overrightarrow{AA_1} = \frac{1}{1+x}\left(x\overrightarrow{AC} + \overrightarrow{AB}\right), \quad \overrightarrow{BB_1} = \frac{1}{1+x}\left(x\overrightarrow{BC} + \overrightarrow{BA}\right)$$

Далее, и из равенства длин чевиан имеем

$$\left(x\overrightarrow{AC} + \overrightarrow{AB}\right)^2 = \left(x\overrightarrow{BC} + \overrightarrow{BA}\right)^2 \Leftrightarrow x^2 AC^2 + AB^2 + 2x\overrightarrow{AC}\cdot\overrightarrow{AB} = x^2 BC^2 + AB^2 + 2x\overrightarrow{BC}\cdot\overrightarrow{BA};$$

$$x(AC^2 - BC^2) + 2(\overrightarrow{AC}\cdot\overrightarrow{AB} + \overrightarrow{BC}\cdot\overrightarrow{AB}) = 0 \Leftrightarrow (\overrightarrow{AC} + \overrightarrow{BC})\cdot(x(\overrightarrow{AC} - \overrightarrow{BC}) + 2\overrightarrow{AB}) = 0;$$

$$2\overrightarrow{MC}\cdot\overrightarrow{AB}(x+2) = 0$$

Значит $x = -2$ или $CM$ перпендикулярна $AB$. Случай $x = -2$ рассмотрен выше (рис. 9). Следовательно, $CM$ является высотой и медианой треугольника $ABC$.

Теорема 2 доказана.

**Замечание 2.1.** Пусть точки $A_1$ и $B_1$ делят стороны $BC$ и $AC$ соответственно в одинаковых отношениях: $BA_1:A_1C=AB_1:B_1C=p:q$, причем $2q + p \neq 0$. Тогда, если $AA_1=BB_1$, то треугольник $ABC$ равнобедренный.

*Доказательство.* Пусть $BA_1:A_1C=AB_1:B_1C =p:q$. Рассмотрим систему материальных точек $\{qA, qB, pC\}$. Так как $2q + p \neq 0$, то существует центр масс этой система – точка $O$. Точки $A_1$ и $B_1$ – центры масс систем материальных точек $\{qA, pC\}$ и $\{qB, pC\}$ соответственно. Центр $O$ лежит на пересечении прямых $AA_1$, $BB_1$ и $CM$, где $M$ – барицентр системы $\{qA, qB\}$. Очевидно, что $CM$ – медиана.

Замечание доказано.

# 5    Доказательство теоремы 3

**Пояснение:** теорему 3 можно доказать двумя способами: алгебраическим и геометрическим. Алгебраическое доказательство заключается в том, что позволяет явно выразить длины чевиан $AA_1$ и $BB_1$, проходящих через точку $T$.

*Алгебраическое доказательство теоремы 3*. Обозначим $\alpha$, $\beta$ и $\gamma$ – углы треугольника $ABC$. Пусть радиус окружности, описанной около треугольника $ABC$ равен 1, а его вершины имеют координаты $A(c_1,0)$, $B(c_2,0)$, $C(0,h)$. Тогда

$$AC = 2\sin\beta,\ BC = 2\sin\alpha,\ AB = 2\sin\gamma,\ c_1 = -2\sin\beta\cos\alpha,$$
$$c_2 = 2\sin\alpha\cos\beta,\ h = 2\sin\alpha\sin\beta,\ \overrightarrow{BC} = (h,c_2)^{\mathrm{T}},\ \overrightarrow{AC} = (h,c_1)^{\mathrm{T}}$$

где операция $(*)^{\mathrm{T}}$ означает транспонирование.

Уравнение окружности, симметричной относительно $AB$ окружности описанной около треугольника $ABC$ запишем в параметрическом виде:

$$\begin{cases} x = \sin t + \dfrac{c_1+c_2}{2},\ 0 \le t \le 2\pi \\ y = \cos t - \cos\gamma \end{cases}$$

Теперь, из векторных уравнений прямых ($AT$) и ($BC$) имеем

$$\overrightarrow{AT}\tau + A = \overrightarrow{BC}\tau_1 + B,\ \overrightarrow{AT}\tau = \overrightarrow{BC}\tau_1 + \overrightarrow{AB}$$

Умножим обе части этого равенства на вектор $(h,c_2)^{\mathrm{T}}$ – ортогональный вектору $\overrightarrow{BC}$:ё

$$\tau\begin{pmatrix} \sin t + \dfrac{c_1+c_2}{2} - c_1 \\ \cos t - \cos\gamma \end{pmatrix}\begin{pmatrix} h \\ c_2 \end{pmatrix} = \begin{pmatrix} 2\sin\gamma \\ 0 \end{pmatrix}\begin{pmatrix} h \\ c_2 \end{pmatrix},\ \tau(h\sin t + h\sin\gamma + c_2\cos t - c_2\cos\gamma) = 2h\sin\gamma$$

$$\tau\cdot 2\sin\alpha(\sin\beta\sin t + \cos\beta\cos t + \sin\beta\sin\gamma - \cos\beta\cos\gamma) = 4\sin\alpha\sin\beta\sin\gamma$$

$$\tau\cdot(\cos(t-\beta) - \cos(\gamma+\beta)) = 2\sin\beta\sin\gamma,\ \tau\cdot(\cos(t-\beta) + \cos\alpha) = 2\sin\beta\sin\gamma,$$

$$\tau\cos\dfrac{\alpha-t+\beta}{2}\cos\dfrac{t+\alpha-\beta}{2} = \sin\beta\sin\gamma,\ \tau = \dfrac{\sin\beta\sin\gamma}{\sin\dfrac{\gamma+t}{2}\cos\dfrac{t+\alpha-\beta}{2}}$$

Отсюда

$$AA_1^2 = \tau^2 AT^2 = \dfrac{\sin^2\beta\sin^2\gamma}{\sin^2\dfrac{\gamma+t}{2}\cos^2\dfrac{t+\alpha-\beta}{2}}((\sin t + \sin\gamma)^2 + (\cos t - \cos\gamma)^2) =$$

$$= \dfrac{\sin^2\beta\sin^2\gamma}{\sin^2\dfrac{\gamma+t}{2}\cos^2\dfrac{t+\alpha-\beta}{2}}(2 - 2\cos(t+\gamma)) = \dfrac{4\sin^2\beta\sin^2\gamma\sin^2\dfrac{\gamma+t}{2}}{\sin^2\dfrac{\gamma+t}{2}\cos^2\dfrac{t+\alpha-\beta}{2}} = \dfrac{4\sin^2\beta\sin^2\gamma}{\cos^2\dfrac{t+\alpha-\beta}{2}}$$

Аналогично находим $BB_1$ и окончательно получаем:

$$BB_1 = \frac{4\sin^2\alpha \sin^2\gamma}{\cos^2\frac{t+\alpha-\beta}{2}} \qquad \text{и} \qquad AA_1 = \frac{4\sin^2\beta \sin^2\gamma}{\cos^2\frac{t+\alpha-\beta}{2}}$$

Отсюда следует, что если $AA_1 = BB_1$, то $\alpha = \beta$.

Теорема 3 доказана.

Для геометрического доказательства этой теоремы понадобятся следующая лемма.

**Лемма 3.1.** *Пусть треугольники ABC и $A_1B_1C_1$ имеют равные острые углы $\angle C = \angle C_1$, равные стороны $AB = A_1B_1$ (или $CB = C_1B_1$) и $\angle A + \angle A_1 = \pi$. Тогда $CB = C_1B_1$ (или $AB = A_1B_1$).*
Доказательство этой леммы следует из рисунка 11.

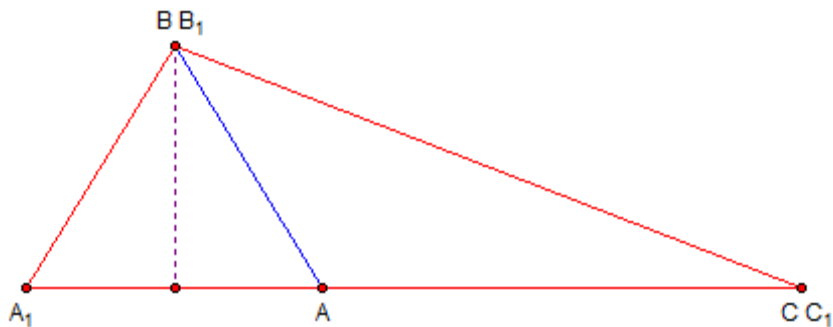
Рис. 11

*Геометрическое доказательство теоремы 3.* Пусть $T$ – точка пересечения равных чевиан треугольника, лежащая на окружности $C_{ABC}$. Сначала заметим, что по условию теоремы точка $T$ не является симметричной вершине $C$ относительно середины отрезка $AB$ (в противном случае $(AT)$ и $(BT)$ будут параллельными сторонам $BC$ и $AC$ соответственно). Случаи $T = A$ и $T = B$ так же исключаются из рассмотрения по понятным причинам. Из всех возможных конфигураций рассмотрим три, изображенных на рисунке 12, где $A$, $B$, $C$ – вершины треугольника $ABC$, $T$ – точка на окружности, $AA_1$, $BB_1$ – чевианы, проходящие через $T$, $C_1$ – точка симметричная точке $C$ относительно $(AB)$.

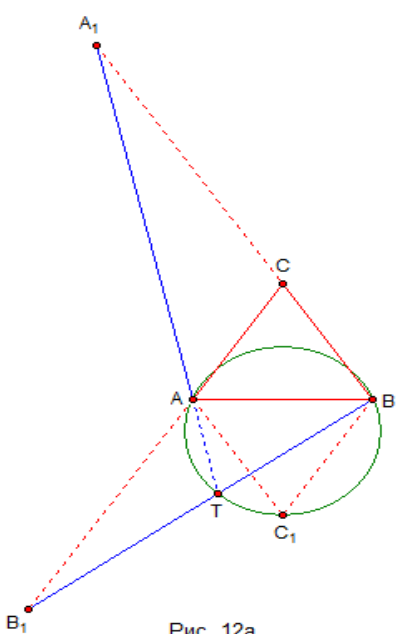
Рис. 12а

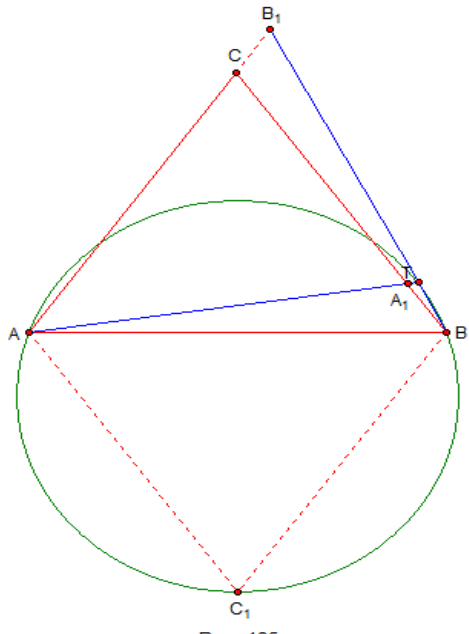
Рис. 12б

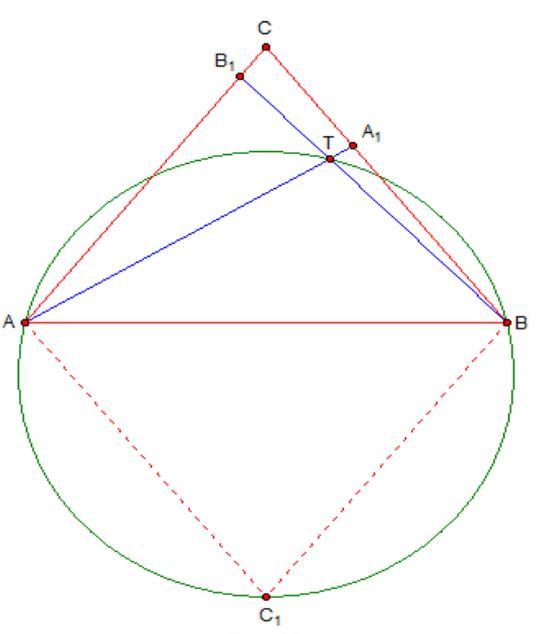
Рис. 12в

Для случая а), когда точка *T* лежит вне треугольника и под стороной *AB*, рассмотрим треугольники $AA_1C$ и $AB_1T$. Углы $B_1AT$ и $A_1AC$ равны как вертикальные. Так как $\angle ATB$ и $\angle ACB$ вписанные и опираются на дугу *AB*, то $\angle ATB_1 = \angle ACA_1$, значит $\angle AB_1T = \angle AA_1C$. Отсюда следует, что треугольники $AA_1C$ и $BB_1C$ имеют равные углы ($\angle AB_1T = \angle AA_1C$), равные стороны ($BB_1=AA_1$) и $\angle ACA_1+\angle BCB_1= \pi$. По лемме 3.1 стороны *BC* и *AC* равны.

В случае б) рассмотрим треугольники $BA_1T$ и $BB_1C$. Они имеют общий угол $A_1BT$ и равные углы $BTA_1$ и $BCB_1$, значит $\angle CB_1B = \angle TA_1B = \angle AA_1C$. Следовательно, треугольники $ACA_1$ и $BCB_1$ имеют равные углы ($\angle AA_1C = \angle BB_1C$), равные стороны $AA_1$ и $BB_1$ и $\angle ACA_1+\angle BCB_1= \pi$. По лемме 3.1 стороны *BC* и *AC* равны.

В случае с) заметим, что сумма углов $B_1TA_1$ и $B_1CA_1$ равна $\pi$. Значит треугольники $ACA_1$ и $BCB_1$ имеют общий угол *C*, равные стороны $AA_1$ и $BB_1$ и $\angle AA_1C+\angle BB_1C= \pi$. По лемме 3.1 стороны *BC* и *AC* равны.

Теорема 3 доказана.

# 6 Доказательство теоремы 4

*Доказательство теоремы 4.* Без потери общности можно считать, что CH = 1. Тогда $AH = ctg\alpha$, $BH = ctg\beta$ и $y = HO$.

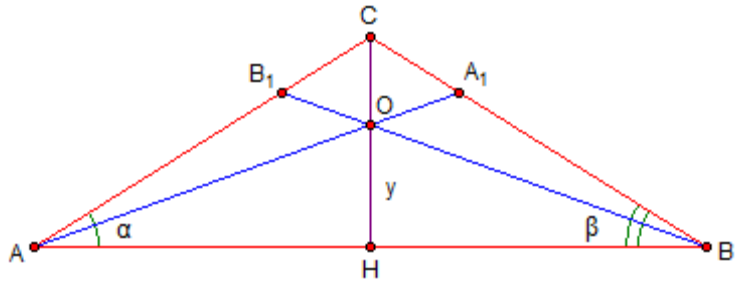

Далее, из теоремы синусов для треугольника $ABB_1$ имеем

$$\frac{BB_1}{\sin BAC} = \frac{AB}{\sin AB_1B} \Rightarrow \frac{BB_1}{\sin \alpha} = \frac{AB}{\sin\left(\alpha + arctg(y \cdot tg\beta)\right)} \Rightarrow BB_1 = \frac{AB\sin \alpha}{\sin\left(\alpha + arctg(y \cdot tg\beta)\right)}$$

Аналогично, в треугольнике $BAA_1$

$$AA_1 = \frac{AB\sin \beta}{\sin\left(\beta + arctg(y \cdot tg\alpha)\right)}$$

Из равенства чевиан получаем:

$$\frac{\sin\left(\beta + arctg(y \cdot tg\alpha)\right)}{\sin \beta} = \frac{\sin\left(\alpha + arctg(y \cdot tg\beta)\right)}{\sin \alpha}$$

После некоторых преобразований приходим к уравнению:

$$\frac{ctg\alpha + y \cdot ctg\beta}{\sqrt{ctg^2\alpha + y^2}} = \frac{ctg\beta + y \cdot ctg\alpha}{\sqrt{ctg^2\beta + y^2}}$$

Возведя в квадрат обе части этого уравнения получаем:

$$(ctg\alpha + y \cdot ctg\beta)^2(ctg^2\beta + y^2) = (ctg\beta + y \cdot ctg\alpha)^2(ctg^2\alpha + y^2);$$
$$y(ctg^2\alpha - ctg^2\beta)(y^3 + y(ctg^2\alpha + ctg^2\beta - 1) + 2ctg\alpha \cdot ctg\beta) = 0$$

Из определения чевиан $O \neq H$, значит $ctg^2\alpha - ctg^2\beta = 0$, отсюда $AH = BH$ и, следовательно, $CH$ является не только высотой, но и медианой треугольника $ABC$, поэтому, треугольник $ABC$ равнобедренный, если только последняя скобка не равна нулю.

Разберемся в тех случаях, когда чевианы равны и пересекаются на высоте неравнобедренного треугольника, тогда:

$$y^3 + y(ctg^2\alpha + ctg^2\beta - 1) + 2ctg\alpha \cdot ctg\beta = 0 \qquad (2)$$

Введем следующие обозначения:

$$\begin{cases} ctg^2\alpha + ctg^2\beta = v \\ ctg\alpha \cdot ctg\beta = u \end{cases} \qquad (3)$$

Тогда уравнение (2) примет вид

$$y^3 + y(v-1) + 2u = 0$$

Известно [2], что если дискриминант $D = -4 \cdot b^3 \cdot d + b^2 \cdot c^2 - 4 \cdot a \cdot c^3 + 18 \cdot a \cdot b \cdot c \cdot d - 27 \cdot a^2 \cdot d^2$ кубического уравнения: $ay^3 + by^2 + cy + d = 0$ больше нуля, то уравнение имеет три различных действительных корня; если же $D = 0$, то оно имеет три действительных корня, причем два из них равны. Если же $D < 0$, то уравнение имеет один корень.

Рассмотри первый случай, когда $D < 0$:

$$D = -4(27u^2 + (v-1)^3) < 0 \Leftrightarrow v > 1 - 3u^{\frac{2}{3}}$$

Учитывая (3), получаем следующее условие:

$$ctg^2\alpha + ctg^2\beta > 1 - 3(ctg\alpha \cdot ctg\beta)^{\frac{2}{3}}$$

Если же неравенство не выполняется, то многочлен

$$p(y) = y^3 + y(ctg^2\alpha + ctg^2\beta - 1) + 2ctg\alpha \cdot ctg\beta$$

имеет три действительных корня $y_0$, $y_1$, $y_2$. Далее производная $p(y)$ должна иметь два корня, не сложно видеть, что они разных знаков, следовательно корни $y_{1,2} > 0$, так как по теореме Виета произведение корней $y_0 y_1 y_2 = -2c_1^2 c_0^2 < 0$. Кроме того, $y_{1,2} < h$. В самом деле, если $y \geq h$, то

$$y^3 + y(c_1^2 + c_0^2 - c_1^2 c_0^2) + 2c_1^2 c_0^2 > y^3 - yh^2 = y(y^2 - h^2) \geq 0$$

Таким образом, мы доказали, что если углы $\alpha$ и $\beta$ не удовлетворяют неравенству:

$$ctg^2\alpha + ctg^2\beta > 1 - 3(ctg\alpha \cdot ctg\beta)^{\frac{2}{3}}$$

то на высоте $CH$ треугольника с такими углами при основании $AB$ существуют две (или три) точки, через которые проходят равные чевианы $AA_1 = BB_1$ и $AA_2 = BB_2$. Треугольник с углами $\alpha=85^0$, $\beta=60^0$ изображён на рисунке 14.

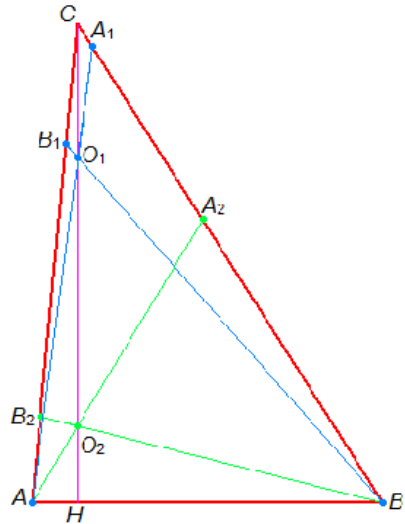

Рис. 14

Учитывая то, что точка $O$ может лежать под основанием треугольника, запишем $y$ в виде:

$$y = \frac{\overrightarrow{HO} \cdot \overrightarrow{HC}}{|HC|}$$

Теорема 4 доказана.

**Следствие 4.1.** На прямой, содержащей высоту $CH$, любого треугольника $ABC$ с острыми углами $\alpha$ и $\beta$ при основании $AB$ существует точка $O$ такая, что проходящие через неё чевианы $AA_1$ и $BB_1$ равны, причём $y = \dfrac{\overrightarrow{HO} \cdot \overrightarrow{HC}}{|HC|}$ удовлетворяет уравнению:

$$y^3 + y(ctg^2\alpha + ctg^2\beta - 1) + 2ctg\alpha \cdot ctg\beta = 0$$

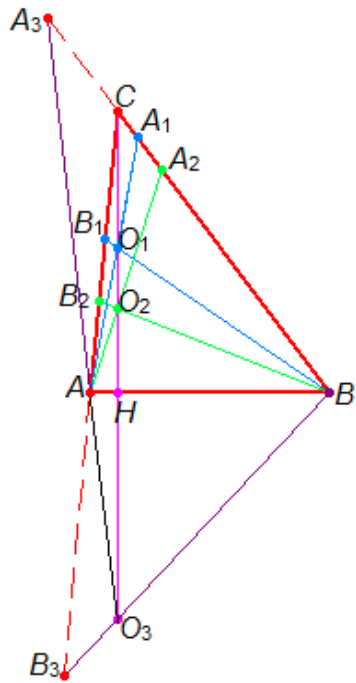

Рис. к следствию 4.1

Если в треугольнике *ABC* угол *A* тупой, то на высоте *CH* существует точка *O*, через которую проходят две равные чевианы $AA_1$ и $BB_1$. В самом деле при движении точки *O* от вершины *C* к основанию, чевиана $BB_1$ будет меняться от *BC* до *BA*, а чевиана $AA_1$ – от *AC* < *BC* до бесконечности (приближаясь к прямой параллельной *BC*), поэтому найдется точка *O*, доказывающее это утверждение (рис. 16).

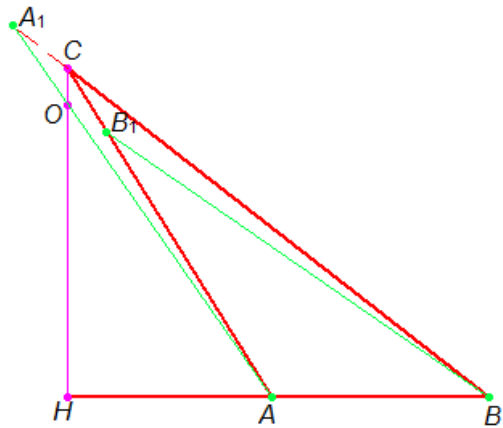

Рис. 16

Пользуясь такими же рассуждениями, как и при доказательстве теоремы 1 легко получить следующее утверждение.

**Следствие 4.2.** На прямой, содержащей высоту *CH*, любого треугольника *ABC* с тупым углом A, существует точка *O* такая, что проходящие через нее чевианы равны и $y = \dfrac{\overrightarrow{HO} \cdot \overrightarrow{HC}}{|HC|}$ удовлетворяет уравнению:

$$y^3 + y(ctg^2\alpha + ctg^2\beta - 1) + 2ctg\alpha \cdot ctg\beta = 0$$

В отличие от уравнения следствия 4.1 это уравнение всегда имеет один положительный корень и может иметь еще два отрицательных (рис. 17).

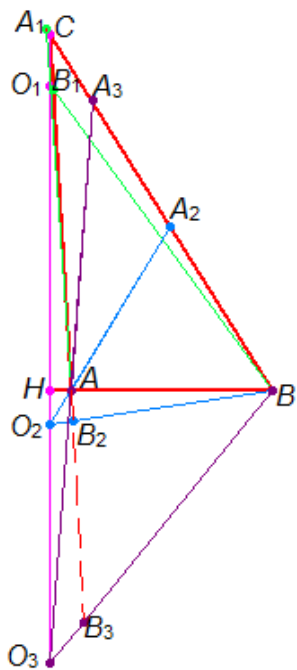

Рис. 17

## 7     Доказательство теоремы 5

**Пояснение:** очевидно, что если чевианы $AA_1$ и $BB_1$ являются 0-трисами углов $BAC$ и $ABC$, то они совпадают со стороной $AB$; если же они являются 1-трисами, то совпадут со сторонами $AC$ и $BC$ соответственно. Если $0 < k < 1$, то точка пересечения таких чевиан лежит внутри треугольника $ABC$ и если $k < 0$ или $k > 1$ точка пересечения – лежит вне $ABC$.

*Доказательство теоремы 5.* Обозначим $\angle A = \alpha$, $\angle B = \beta$ и $\angle C = \gamma$. По теореме синусов для треугольников $ABB_1$ и $BAA_1$ имеем:

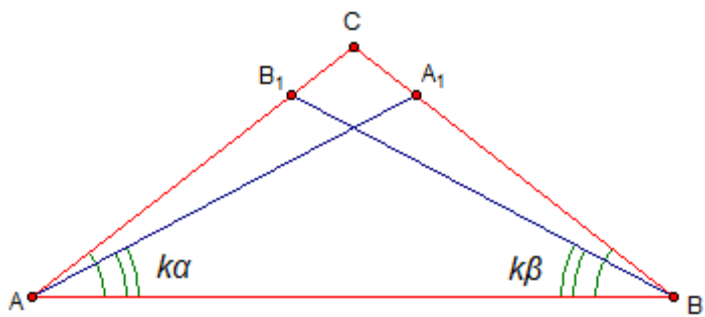

$$\frac{AB}{\sin(k\alpha+\beta)}=\frac{AA_1}{\sin\beta},$$

$$\frac{AB}{\sin(k\beta+\alpha)}=\frac{BB_1}{\sin\alpha}$$

Так как $AA_1=BB_1$, то отсюда имеем

$$\sin(k\alpha+\beta)\sin\alpha = \sin(k\beta+\alpha)\sin\beta,$$
$$\cos((k-1)\alpha+\beta)-\cos((k+1)\alpha+\beta)=\cos((k-1)\beta+\alpha)-\cos((k+1)\beta+\alpha),$$
$$\cos((k-1)\alpha+\beta)-\cos((k-1)\beta+\alpha)=\cos((k+1)\alpha+\beta)-\cos((k+1)\beta+\alpha),$$
$$\sin\frac{(\alpha+\beta)k}{2}\sin\frac{(\alpha-\beta)(k-2)}{2}=\sin\frac{(\alpha+\beta)(k+2)}{2}\sin\frac{(\alpha-\beta)k}{2}$$

(Это равенство будет выполняться при $\alpha = \beta$). Будем считать, что $\alpha > \beta$. Обозначим $u=\dfrac{\alpha+\beta}{2}$ и $v=\dfrac{\alpha-\beta}{2}$. Тогда

$$\begin{cases}\dfrac{\sin(k+2)u}{\sin ku}=\dfrac{\sin(k-2)v}{\sin kv}\\ 0<v<u<\dfrac{\pi}{2}\end{cases}$$

Положим $t=ku$ и, $\tau=kv$ и $\gamma=\dfrac{2}{k}$. Тогда наше равенство примет вид:

$$\frac{\sin(\gamma+1)t}{\sin t}=-\frac{\sin(\gamma-1)\tau}{\sin\tau}$$

Очевидно, что $0<\tau<t<\dfrac{\pi}{\gamma}$ по определению, а $\gamma \geq 2$ – по условию теоремы.

Рассмотрим функции

$$f(t)=\frac{\sin(\gamma+1)t}{\sin t} \quad \text{и} \quad g(\tau)=-\frac{\sin(\gamma-1)\tau}{\sin\tau} \qquad (4)$$

где $0<\tau<t<\dfrac{\pi}{\gamma}$. Найдем общее значение, которое принимают эти функции: $f\left(\dfrac{\pi}{\gamma}\right)=g\left(\dfrac{\pi}{\gamma}\right)=-1$.

Причем, на промежутке $0<\tau<t<\dfrac{\pi}{\gamma}$ функция $f(t)=\dfrac{\sin(\gamma+1)t}{\sin t}>-1$, а $g(\tau)=-\dfrac{\sin(\gamma-1)\tau}{\sin\tau}<-1$.

Следовательно, при $0<\tau<t<\dfrac{\pi}{\gamma}$ функции $f(t)$ и $g(\tau)$ не принимают одинаковых значений. Поэтому, $\alpha=\beta$ и треугольник равнобедренный.

Это и доказывает теорему 5.

**Замечание 5.1.** Если $k$-трисы ($AA_1$) и ($BB_1$) не пересекают соответственно стороны $BC$ и $AC$ треугольника $ABC$, то треугольник $ABC$ равнобедренный.

Действительно, так как $AA_1 \parallel BC$ и $BB_1 \parallel AC$ и $(AA_1)$ и $(BB_1)$ – $k$-трисы, то

$$\begin{cases} k\alpha + \beta = \pi \\ k\beta + \alpha = \pi \end{cases} \Rightarrow \alpha = \beta$$

Значит, треугольник $ABC$ – равнобедренный.

При $k > 1$ легко показать, что в любом прямоугольном треугольнике ($\angle C = 90°$) 2-трисы равны и равны стороне $AB$. Однако, имеет место следующие утверждения.

**Утверждение 5.2**. Не прямоугольный треугольник с равными 2-трисами ($AA_1=BB_1$) равнобедренный.
Доказательство утверждения следует из равенств (1) при $k = 2$.

**Утверждение 5.3**. Для любого $k > 1$ существуют неравнобедренные треугольники с равными $k$-трисами.
Доказательство утверждения следует из свойств функций (4).

## 8   Доказательство теоремы 6

*Доказательство теоремы 6.* Пусть $AA_1= AA_2= BB_1= BB_2= CC_1= CC_2 = l$. Обозначим также $h_a$, $h_b$ и $h_c$, – высоты треугольника $ABC$, опущенные на стороны $BC$, $AC$ и $AB$ – соответственно, а $H_a$, $H_b$ и $H_c$ – основания этих высот.

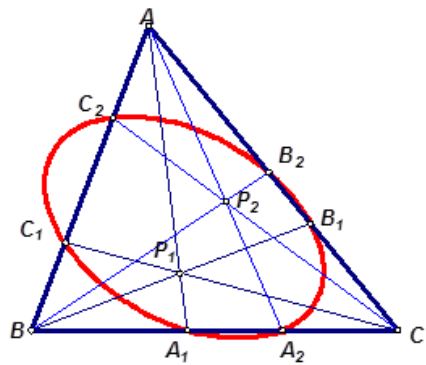

Из условия теоремы следует, что треугольники $AA_1A_2$, $BB_1B_2$ и $CC_1C_2$ равнобедренные, поэтому, например, в треугольнике $CC_1C_2$

$$C_1H_c = C_2H_c$$

отсюда

$$BC_1 \cdot BC_2 = (BH_c - C_1H_c) \cdot (BH_c + C_1H_c) = BH_c^2 - C_1H_c^2 = BH_c^2 - (l^2 - h_c^2) =$$
$$= BH_c^2 + h_c^2 - l^2 = BC^2 - l^2$$

Аналогично, получаем:

$$CA_1 \cdot CA_2 = AC^2 - l^2; \ AC_1 \cdot AC_2 = AC^2 - l^2; \ BA_1 \cdot BA_2 = AB^2 - l^2;$$
$$AB_1 \cdot AB_2 = AB^2 - l^2; CB_1 \cdot CB_2 = BC^2 - l^2$$

Проверим условие Карно:

$$\left(\frac{BA_1}{CA_1} \cdot \frac{BA_2}{CA_2}\right) \cdot \left(\frac{CB_1}{AB_1} \cdot \frac{CB_2}{AB_2}\right) \cdot \left(\frac{AC_1}{BC_1} \cdot \frac{AC_2}{BC_2}\right) = \frac{AB^2 - l^2}{AC^2 - l^2} \cdot \frac{BC^2 - l^2}{AB^2 - l^2} \cdot \frac{AC^2 - l^2}{BC^2 - l^2} = 1$$

Теорема доказана.

# 9 Доказательство теоремы 7

Найдем геометрическое место точек пересечения равных чевиан произвольного треугольника *ABC*, проведенных из вершин *A* и *B*. Сразу же отметим, что все точки, лежащие на основание *AB*, принадлежат этому множеству.

Для треугольника *ABC* введем систему координат с центром в точке *M*, где *CM* – медиана треугольника *ABC*, а ось абсцисс содержит основание *AB*. Тогда, из теоремы синусов, координаты вершин треугольника *ABC* имеют вид:

$$A(-c, 0), \ B(c, 0),$$
$$C\left(c - \frac{2c \sin\alpha \cos\beta}{\sin\gamma}, \frac{2c \sin\alpha \sin\beta}{\sin\gamma}\right) = C\left(\frac{c \sin(\beta - \alpha)}{\sin(\alpha + \beta)}, \frac{2c \sin\alpha \sin\beta}{\sin(\alpha + \beta)}\right)$$

где *c* > 0 – некоторое число, а *α*, *β*, *γ* – углы треугольника *ABC*.

Если $AA_1$ и $BB_1$ – равные чевианы треугольника *ABC*, а $(x_0, y_0)$ – точка их пересечения, то из параметрических уравнений прямых ($AA_1$) и (*BC*) имеем

$$\begin{pmatrix} x_0 + c \\ y_0 \end{pmatrix} \tau + \begin{pmatrix} -c \\ 0 \end{pmatrix} = \begin{pmatrix} -\cos\beta \\ \sin\beta \end{pmatrix} t + \begin{pmatrix} c \\ 0 \end{pmatrix} \quad \left[\text{умножим на } \begin{pmatrix} \sin\beta \\ \cos\beta \end{pmatrix}\right]$$

$$((x_0 + c)\sin\beta + y_0 \cos\beta)\tau = 2c \sin\beta \Rightarrow \tau = \frac{2c \sin\beta}{(x_0 + c)\sin\beta + y_0 \cos\beta} \Rightarrow$$

$$|AA_1|^2 = \left|\begin{pmatrix} x_0 + c \\ y_0 \end{pmatrix} \frac{2c \sin\beta}{(x_0 + c)\sin\beta + y_0 \cos\beta}\right|^2 = \frac{4c^2 \sin^2\beta((x_0 + c)^2 + y_0^2)}{(y_0 \cos\beta + (x_0 + c)\sin\beta)^2}$$

Аналогично найдем квадрат длины $BB_1$:

$$|BB_1|^2 = \frac{4c^2 \sin^2\alpha((x_0 - c)^2 + y_0^2)}{(y_0 \cos\alpha - (x_0 - c)\sin\alpha)^2}$$

Отсюда, учитывая равенство этих чевиан, приходим к уравнению:

$$\frac{\sin^2 \beta((x_0+c)^2 + y_0^2)}{(y_0\cos\beta + (x_0+c)\sin\beta)^2} = \frac{\sin^2 \alpha((x_0-c)^2 + y_0^2)}{(y_0\cos\alpha - (x_0-c)\sin\alpha)^2}$$

После несложных преобразований и, полагая, для удобства, $x_0=x$, $y_0=y$ и $2c=1$, приходим к уравнению:

$$y((x^2+y^2)(y\sin(\beta-\alpha)\sin(\beta+\alpha) - 2x\sin\alpha\sin\beta) +$$
$$+ (y^2-x^2)\sin\alpha\sin\beta\sin(\beta-\alpha) + xy(\sin^2\alpha\cos^2\beta + \sin^2\beta\cos^2\alpha -$$
$$- 2\sin^2\beta\sin^2\alpha) + \frac{1}{2}x\sin\alpha\sin\beta\sin(\beta+\alpha) + \frac{1}{4}y\sin(\beta-\alpha)\sin(\beta+\alpha) +$$
$$+ \frac{1}{4}\sin\alpha\sin\beta\sin(\beta-\alpha)) = 0$$

Отсюда получаем первое (тривиальное) решение: $y=0$ и кубическую кривую, определяемую уравнением:

$$(x^2+y^2)\sin(\beta+\alpha)(y\sin(\beta-\alpha) - 2x\sin\alpha\sin\beta) +$$
$$+ (y^2-x^2)\sin\alpha\sin\beta\sin(\beta-\alpha) + xy(\sin^2\alpha\cos^2\beta + \sin^2\beta\cos^2\alpha -$$
$$- 2\sin^2\beta\sin^2\alpha) + \frac{1}{2}x\sin\alpha\sin\beta\sin(\beta+\alpha) + \frac{1}{4}y\sin(\beta-\alpha)\sin(\beta+\alpha) + \quad (5)$$
$$+ \frac{1}{4}\sin\alpha\sin\beta\sin(\beta-\alpha) = 0$$

Как известно, кривые третьего порядка имеют довольно сложные классификации [7], [8]. Отметим некоторые ее свойства.
В [7, стр. 44] показано, что если кубическая кривая имеет вид:

$$Ax^3 + 3Bx^2y + 3Cy^2x + Dy^3 + 3Ex^2 + 6Fxy + 3Gy^2 + 3Hx + 3Ky + L = 0$$

а $y=kx+b$ – уравнение ее асимптоты, то угловой коэффициент $k$ и начальная ордината $b$ определяются равенствами:

$$A + 3Bk + 3Ck^2 + Dk^3 = 0, \quad (B + 2Ck + Dk^2)b = -(E + 2Fk + Gk^2)$$

при этом начальная ордината $b$ может не существовать. Но, в нашем случае $k$ и $b$ находятся легко и определяются равенствами:

$$(k^2+1)(k\sin(\beta-\alpha) + 2\sin\alpha\sin\beta) = 0, \quad k = \frac{2\sin\alpha\sin\beta}{\sin(\beta-\alpha)}$$
$$b = -\frac{2\sin\alpha\sin\beta\sin(\beta+\alpha)}{\sin^2(\beta-\alpha) + 4\sin^2\alpha\sin^2\beta}$$

Таким образом, наша кривая имеет одну асимптоту:

$$y = \frac{2\sin\alpha\sin\beta}{\sin(\beta-\alpha)} x - \frac{2\sin\alpha\sin\beta\sin(\beta+\alpha)}{\sin^2(\beta-\alpha) + 4\sin^2\alpha\sin^2\beta} \quad (6)$$

Следуя классификации Ньютона [7, стр. 47] эта кривая относится ко второй группе, имеет одну асимптоту и одну бесконечную ветвь прямолинейного типа (прямолинейной называется гиперболическая ветвь, вытянутая вдоль прямой, являющейся асимптотой).
Таким образом, мы приходим к следующему утверждению

Геометрическое место точек пересечения пар равных чевиан $AA_1$ и $BB_1$ треугольника $ABC$ есть объединение отрезка $AB$ и кубической кривой. Эта кривая обладает следующими свойствами:
 а) пересекает ($AB$) ровно в трех точках: $A$, $B$ и $E$, причем Е симметрична относительно $M$ – середины $AB$ точке $H$ – основанию высоты $CH$ треугольника $ABC$;
 б) пересекает серединный перпендикуляр стороны $AB$ только в одной точке $N$, лежащей ниже основания $AB$ и равным высоте $CH$ треугольника $ABC$ ($MN=CH$);
 в) ее асимптота параллельна медиане $CM$ треугольника $ABC$;
 г) точка $D$ самопересечения (узел) гиперболы симметрична вершине $C$ относительно $M$ (середины $AB$).

*Доказательства*:
*а*) Заметим, что $A(-\frac{1}{2},0)$, $B(\frac{1}{2},0)$, $C\left(\frac{\sin(\beta-\alpha)}{2\sin(\alpha+\beta)}, \frac{\sin\alpha\sin\beta}{\sin(\alpha+\beta)}\right)$ теперь, подставляя $y = 0$ в равнение (5), получим утверждение *а*):

$$-2x^3\sin(\beta+\alpha)\sin\alpha\sin\beta - x^2\sin\alpha\sin\beta\sin(\beta-\alpha) + \frac{1}{2}x\sin\alpha\sin\beta\sin(\beta+\alpha) +$$
$$+ \frac{1}{4}\sin\alpha\sin\beta\sin(\beta-\alpha) = 0$$
$$x^3\sin(\beta+\alpha) + \frac{1}{2}x^2\sin(\beta-\alpha) - \frac{1}{4}x\sin(\beta+\alpha) - \frac{1}{8}\sin(\beta-\alpha) = 0$$
$$x\sin(\beta+\alpha)\left(x^2 - \frac{1}{4}\right) + \frac{1}{2}\sin(\beta-\alpha)\left(x^2 - \frac{1}{4}\right) = 0; \quad x_1 = -\frac{1}{2}; \; x_2 = \frac{1}{2}; \; x_3 = -\frac{\sin(\beta-\alpha)}{2\sin(\beta+\alpha)}$$

Аналогично докажем б), подставляя в (5) $x=0$:

$$y^3\sin(\beta+\alpha)\sin(\beta-\alpha) + y^2\sin\alpha\sin\beta\sin(\beta-\alpha) + \frac{1}{4}y\sin(\beta-\alpha)\sin(\beta+\alpha) +$$
$$+ \frac{1}{4}\sin\alpha\sin\beta\sin(\beta-\alpha) = 0$$
$$y^3\sin(\beta+\alpha) + y^2\sin\alpha\sin\beta + \frac{1}{4}y\sin(\beta+\alpha) + \frac{1}{4}\sin\alpha\sin\beta = 0$$
$$\left(y + \frac{\sin\alpha\sin\beta}{\sin(\beta+\alpha)}\right)\left(y^2 + \frac{1}{4}\right) = 0 \Rightarrow y_N = -\frac{\sin\alpha\sin\beta}{\sin(\beta+\alpha)} = -y_C$$

Для доказательства в) найдем котангенс угла *CMB*:

$$\operatorname{ctg} CMB = \frac{\dfrac{1}{2} - \dfrac{\sin\alpha\cos\beta}{\sin\gamma}}{\dfrac{\sin\alpha\sin\beta}{\sin\gamma}} = \frac{\sin(\alpha+\beta) - 2\sin\alpha\cos\beta}{2\sin\alpha\cos\beta} = \frac{\sin(\beta-\alpha)}{2\sin\alpha\cos\beta}$$

Отсюда и (6) получаем утверждение г).

Пункт г) достаточно очевиден, т.к. *ACBD* – параллелограмм. То есть, при приближении к точке D, чевианы становятся параллельными сторонам треугольника (их длины бесконечно возрастают). Однако приведем более строгие рассуждения, полезные для построения нашей кривой.

Для построения полученной гиперболы найдем ее параметрическое представление. В качестве параметра будем использовать величину *l* – длины чевиан. Очевидно, что $l \geq \max\{h_a, h_b\}$, где $h_a$ и $h_b$ – высоты треугольника *ABC*, опущенные из вершин *A* и *B* соответственно.

Пусть $AA_1$, $AA_2$, $BB_1$ и $BB_2$ – чевианы треугольника *ABC*, длины которых равны *l*. Тогда, легко видеть, что т.к. *AB*=1, то

$$AB_i = \cos\alpha \pm \sqrt{l^2 - \sin^2\alpha}, \quad BA_i = \cos\beta \pm \sqrt{l^2 - \sin^2\beta}, \quad i = 1,2$$

Найдем точки ($x_{i,j}$, $y_{i,j}$) – пересечения чевиан $AA_i$, $BB_j$, используя параметрические уравнения прямых ($AA_i$) и ($BB_j$), $i, j = 1, 2$:

$$\left(AB_i \begin{pmatrix} \cos\alpha \\ \sin\alpha \end{pmatrix} + \begin{pmatrix} -1 \\ 0 \end{pmatrix}\right) t + \begin{pmatrix} 0{,}5 \\ 0 \end{pmatrix} = \left(BA_j \begin{pmatrix} -\cos\beta \\ \sin\beta \end{pmatrix} + \begin{pmatrix} 1 \\ 0 \end{pmatrix}\right) \tau + \begin{pmatrix} -0{,}5 \\ 0 \end{pmatrix}$$

$$\begin{pmatrix} AB_i \cos\alpha - 1 \\ AB_i \sin\alpha \end{pmatrix} t + \begin{pmatrix} 1 \\ 0 \end{pmatrix} = \begin{pmatrix} -BA_j \cos\beta + 1 \\ BA_j \sin\beta \end{pmatrix} \tau$$

$$\left(\begin{pmatrix} AB_i \cos\alpha - 1 \\ AB_i \sin\alpha \end{pmatrix} t + \begin{pmatrix} 1 \\ 0 \end{pmatrix}\right) \cdot \begin{pmatrix} BA_j \sin\beta \\ BA_j \cos\beta - 1 \end{pmatrix} = 0$$

$$t = \frac{-BA_j \sin\beta}{AB_i BA_j \sin(\alpha+\beta) - (BA_j \sin\beta + AB_i \sin\alpha)}$$

Отсюда получаем параметрическое уравнение нашей кривой
Точнее, мы получили уравнения 4-х кривых, из которых и состоит наша гипербола.

$$\begin{pmatrix} x_{i,j} \\ y_{i,j} \end{pmatrix} = \begin{pmatrix} \dfrac{-BA_j AB_i \sin\beta\cos\alpha + 2cBA_j \sin\beta}{AB_i BA_j \sin(\alpha+\beta) - BA_j \sin\beta - AB_i \sin\alpha} + \dfrac{1}{2} \\ \dfrac{-BA_j AB_i \sin\alpha\sin\beta}{AB_i BA_j \sin(\alpha+\beta) - BA_j \sin\beta - AB_i \sin\alpha} \end{pmatrix}, \quad i, j = 1,2 \quad (7)$$

С помощью этих уравнений легко строить наши гиперболы. Кроме того, из (7) видно, что при неограниченном увеличении длин чевиан *l* полученные 4 кривые пересекаются в узловой точке *D*:

$$\lim_{l\to\infty}\begin{pmatrix}x_{i,j}\\y_{i,j}\end{pmatrix}=\lim_{l\to\infty}\begin{pmatrix}\dfrac{-\sin\beta\cos\alpha+\dfrac{\sin\beta}{AB_i}}{\sin(\alpha+\beta)-\dfrac{\sin\beta}{AB_i}-\dfrac{\sin\alpha}{BA_j}}+\dfrac{1}{2}\\ \dfrac{-\sin\alpha\sin\beta}{\sin(\alpha+\beta)-\dfrac{\sin\beta}{AB_i}-\dfrac{\sin\alpha}{BA_j}}\end{pmatrix}=\begin{pmatrix}\dfrac{-\sin(\beta-\alpha)}{2\sin(\alpha+\beta)}\\ \dfrac{-\sin\beta\cos\alpha}{\sin(\alpha+\beta)}\end{pmatrix}$$

Сравнивая с точкой $C\left(\dfrac{\sin(\beta-\alpha)}{2\sin(\alpha+\beta)},\dfrac{\sin\alpha\sin\beta}{\sin(\alpha+\beta)}\right)$ получаем утверждение д).

На рисунках ниже представлены примеры наших кривых для треугольников с углами $\alpha=20°$, $\beta=40°$ (рис. 18а) и $\alpha=40°$, $\beta=120°$ (рис. 18б).

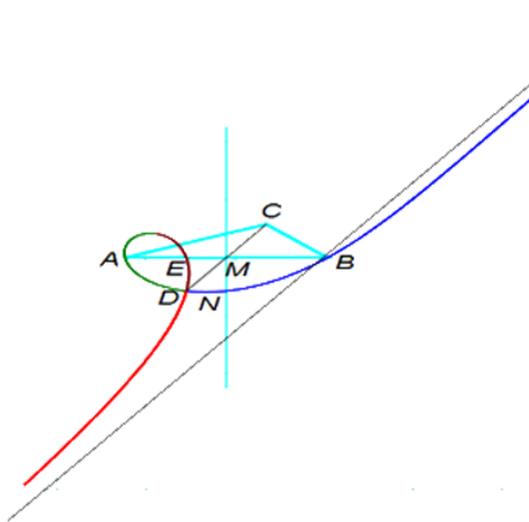
Рис. 18а

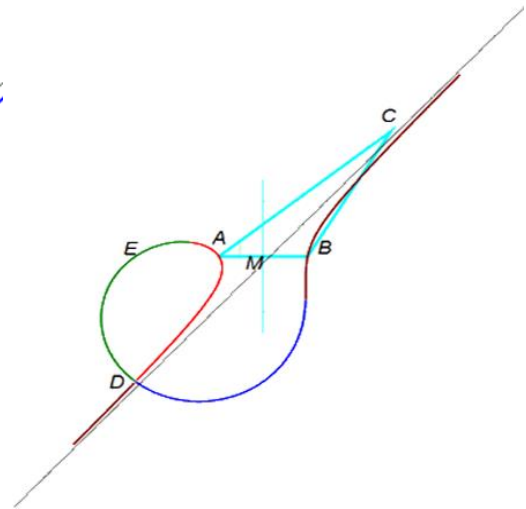
Рис. 18б

В заключении приведем полезное замечание для равнобедренных треугольников

**Замечание 7.1.** Геометрическое место точек пересечения пар равных чевиан $AA_1$ и $BB_1$ равнобедренного треугольника $ABC$ есть объединение отрезка $AB$, оси его симметрии и окружности, симметричной относительно $AB$ окружности, описанной около треугольника $ABC$.

*Доказательство.* Подставляя в (5) $\alpha=\beta$, получим

$$-2x(x^2+y^2)\sin 2\alpha\sin^2\alpha+2xy(\sin^2\alpha\cos^2\alpha-\sin^4\alpha)+\dfrac{1}{2}x\sin^2\alpha\sin 2\alpha=0$$

$$-x(2(x^2+y^2)\sin 2\alpha-2y\cos 2\alpha-\dfrac{1}{2}\sin 2\alpha)=0$$

Отсюда имеем первое решение $x=0$, т.е. ось симметрии треугольника $ABC$. Далее, имеем

$$x^2\sin 2\alpha+y^2\sin 2\alpha-y\cos 2\alpha-\dfrac{1}{4}\sin 2\alpha=0,\ y^2-y\dfrac{\cos 2\alpha}{\sin 2\alpha}+x^2-\dfrac{1}{4}=0$$

$$\left(y+\dfrac{\cos\gamma}{2\sin\gamma}\right)^2+x^2=\left(\dfrac{1}{2\sin\gamma}\right)^2=R^2$$

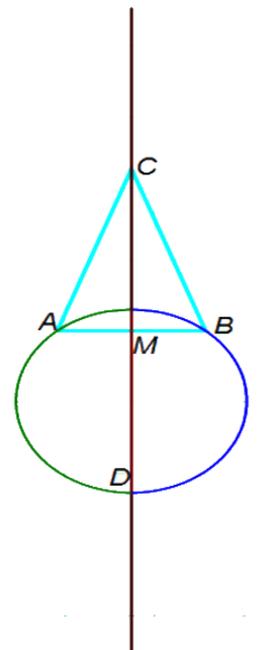
Рис. 19

где *γ=π-2α*, а *R* – радиус, описанной около треугольника *ABC* окружности (рис.19).
   Замечание доказано.

## 10   Переходя к трехмерному случаю

В качестве продолжения исследования можно перейти к пространственному случаю, рассматривая тетраэдр и называя чевианой отрезок, соединяющий его вершину с точкой противоположной грани. Постановка задачи в этом случае может быть сформулирована следующим образом: *выяснить, является ли тетраэдр равногранным (т.е. имеющим равные грани), если все его чевианы равны*.

Для удобства будем пользоваться эквивалентным определением равногранного тетраэдра: доказано, что грани тетраэдра равны тогда и только тогда, когда они равновелики. Начнем рассмотрение поставленной задачи с наиболее простых случаев, когда в качестве чевиан тетраэдра выступают его высоты или биссектрисы трёхгранных углов.

Равногранность тетраэдра с равными высотами сразу же вытекает из формулы для вычисления его объема. Попробуем сделать некоторые шаги к проверке гипотезы о равногранности тетраэдра с равными биссектрисами.

## 11   Пояснение к гипотезе 8

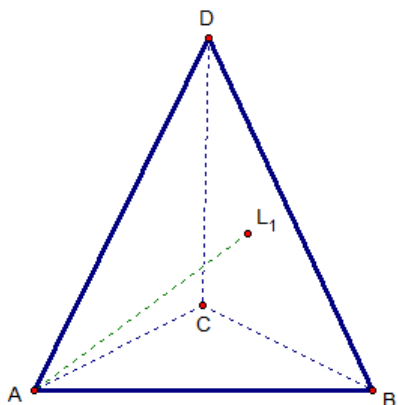

Рис. 20

Биссектрисой трехгранного угла *ABCD* тетраэдра *ABCD* называют отрезок $AL_1$, соединяющий вершину *A* с точкой противоположной грани $L_1$ так, что $\angle L_1AB = \angle L_1AC = \angle L_1AD$ (рис. 20). Точка пересечения всех биссектрис тетраэдра является центром вписанного в него шара.

Длины биссектрис тетраэдра могут быть вычислены по следующей формуле:

$$AL_1 = \frac{(S_{ACD} + S_{ABD} + S_{ABC})(AB^2 \cdot S_{ACD} + AD^2 \cdot S_{ABC} + AC^2 \cdot S_{ABD}) - S_{ACD} \cdot S_{ABC} \cdot DB^2 - S_{ACD} \cdot S_{ABD} \cdot CB^2 - S_{ABD} \cdot S_{ABC} \cdot DA^2}{(S_{ACD} + S_{ABD} + S_{ABC})^2}$$

$$BL_2 = \frac{(S_{ABC} + S_{BCD} + S_{ABD})(DB^2 \cdot S_{ABC} + CB^2 \cdot S_{ABD} + AB^2 \cdot S_{BCD}) - S_{ABC} \cdot S_{ABD} \cdot DA^2 - S_{ABD} \cdot S_{BCD} \cdot AC^2 - S_{ABC} \cdot S_{BCD} \cdot DC^2}{(S_{ACD} + S_{BCD} + S_{ABD})^2}$$

(8)

$$CL_3 = \frac{(S_{ACD} + S_{BCD} + S_{ABC})(CB^2 \cdot S_{ACD} + CD^2 \cdot S_{ABC} + AC^2 \cdot S_{BCD}) - S_{ACD} \cdot S_{ABC} \cdot DB^2 - S_{ABC} \cdot S_{BCD} \cdot DC^2 - S_{ACD} \cdot S_{BCD} \cdot AB^2}{(S_{ACD} + S_{BCD} + S_{ABC})^2}$$

$$DL_4 = \frac{(S_{ACD} + S_{BCD} + S_{ABD})(BD^2 \cdot S_{ACD} + AD^2 \cdot S_{BCD} + DC^2 \cdot S_{ABD}) - S_{ACD} \cdot S_{ABD} \cdot CB^2 - S_{ACD} \cdot S_{BCD} \cdot AB^2 - S_{ABD} \cdot S_{BCD} \cdot AC^2}{(S_{ACD} + S_{BCD} + S_{ABD})^2}$$

Для проверки гипотезы о равногранности тетраэдра с равными биссектрисами был проведён численный эксперимент для тетраэдра $DABC$, у которого углы грани $ABC$ равны $45°$, $60°$, $75°$. Приняв за единицу диаметр описанной около грани $ABC$ окружности, получим: $CB = sin\angle CAB$, $AB = sin\angle ACB$, $AC = sin\angle CBA$. С помощью математического пакета Mathcad была решена система нелинейных относительно $AD = x$, $BD = y$, $CD = z$ уравнений

$$\begin{cases} AL_1 = BL_2 \\ AL_1 = CL_3 \\ AL_1 = DL_4 \end{cases} \text{ с ограничениями: } \begin{cases} DB + DA > CB \\ DB + CB > DA \\ DC + AB > DB \\ DB + AB > AB \\ CB + DA > DB \\ DB + AB > DC \\ DC + DA > AC \\ DC + AC > DA \\ DA + AC > DC \end{cases}$$

В этой системе длины биссектрис $AL_1$, $BL_2$, $CL_3$ и $DL_4$ рассчитывались по формулам (8), в которых площади граней тетраэдра зависят от неизвестных рёбер $x$, $y$ и $z$. Полученное в результате расчёта решение системы $R = \{x, y, z\} = \{8{,}660, 7{,}071, 9{,}659\}$ позволило сделать вывод о равногранности тетраэдра. В частности, площади его граней оказались равными $S_{ACD} = S_{BCD} = S_{ABD} = S_{ABC} = 118{,}301$ из чего следует, что рассмотренный тетраэдр является равногранным. Тем не менее, понятно, что рассмотренный частный случай позволяет лишь выдвинуть гипотезу о равногранности тетраэдра с равными чевианами (или хотя бы биссектрисами), проверка которой должна служить предметом дальнейших изысканий.

Развитием данной работы может стать рассмотрение трёхмерной задачи о формулировке и доказательстве признаков равногранности тетраэдра, связанными с его чевианами.

# Список литературы